\documentclass[12pt]{article}
\title{Extended finite operator calculus as an example of algebraization of
analysis}
\author{A.K.Kwa\'sniewski\\  %%% ENDNAME
\\ Member of the Institute of Combinatorics and its Applications\\
High School of Mathematics and Applied Informatics\\
PL-15-021 Bia{\l}ystok, ul.Kamienna 17, POLAND\\
e-mail: kwandr@gmail.com}

\usepackage{amsmath,amsthm,amsfonts,longtable}

\newtheorem{thm}{Theorem}[section]
\newtheorem{prop}{Proposition}[section]
\newtheorem{note}{Note}[section]
\newtheorem{n}{Notation}[section]
\newtheorem{ex}{Example}[section]
\newtheorem{obs}{Observation}[section]
\newtheorem{com}{Comment}[section]
\newtheorem{rem}{Remark}[section]
\newtheorem{defn}{Definition}[section]
\newtheorem{ident}{Identification}[section]
\numberwithin{equation}{section}

\begin{document}
\maketitle

\begin{abstract}

\noindent ``A Calculus of Sequences'' started in 1936 by Ward constitutes
the general scheme for extensions of classical operator calculus
of Rota - Mullin considered by many afterwards and after Ward.
Because of the notation we shall call the Ward`s calculus of
sequences in its afterwards elaborated form - a $\psi $-calculus.

\noindent The $\psi $-calculus in parts appears to be almost automatic, natural
extension of classical operator calculus of Rota - Mullin or equivalently -
of umbral calculus of Roman and Rota.

\noindent This is a review article based on the $1999-2002$  author`s relevant
contributions.

\end{abstract}

\small{
MCS (2000) : 05A40 , 81S99}

\vspace{0.1cm}

\noindent revised version received 13 July \textbf{2003} and
published in  Bulletin of the Allahabad Mathematical Society  vol. 19, 2005, 1-24.

\vspace{0.1cm}

\noindent affiliated to The Internet Gian-Carlo Polish Seminar:

\noindent http://ii.uwb.edu.pl/akk/sem/sem\_rota.htm

\section{Introduction}

We shall call the Wards calculus of sequences \cite{1} in its
subsequent last century elaborated form - a $\psi $-calculus
because of the notation \cite{2}-\cite{7}. The efficiency of the
Rota oriented language and our notation used has been already
exemplified by simpler proofs of $\psi $-extended counterparts of
all representation independent statements of $\psi $-calculus
\cite{2}. Here these are $\psi $-labelled representations of
Graves-Heisenberg-Weyl (GHW)\cite{13},\cite{1},\cite{16},\cite{15}
algebra of linear operators acting on the algebra $P$ of
polynomials.

As a matter of fact $\psi $-calculus becomes in parts an almost
automatic extension of Rota - Mullin calculus\cite{8} or
equivalently - of the umbral calculus of Roman and Rota
\cite{8,9,10}. The $\psi $-extension relies on the notion of
$\partial _{\psi} $-shift invariance of operators with $\psi
$-derivatives $\partial _{\psi} $ staying for equivalence classes
representatives of special differential operators lowering the
degrees of polynomials by one \cite{6,7,11}. Many of the results
of $\psi $-calculus may be extended to Markowsky $Q$-umbral
calculus \cite{11} where $Q$ stands for arbitrary generalized
difference operator, i.e. the one lowering the degree of any
polynomial by one. $Q$-umbral calculus \cite{11} - as we call it -
includes also those generalized difference operators, which are
not series in $\psi $-derivative $\partial_{\psi}  $ whatever an
admissible $\psi $ sequence would be  (for - "admissible" - see
next section).

The note is at the same time the operator formulation of ``A Calculus of
Sequences'' started in 1936 by Ward \cite{1} with the indication of the role
the $\psi $-representations of Graves-Heisenberg-Weyl (GHW) algebra in
formulation and derivation of principal statements of the $\psi $-extension
of finite operator calculus of Rota.

Restating what was said above we observe that all statements of standard
finite operator calculus of Rota are valid also in the case of $\psi
$-extension under the almost automatic replacement of $\{D,\hat {x}, id \}$
generators of GHW by their $\psi $-representation correspondents
$\{\partial _{\psi},\hat {x}_{\psi}, id\}$ - see definitions 2.1 and
2.5. Naturally any specification of admissible $\psi $ - for example the
famous one defining q-calculus - has its own characteristic properties not
pertaining to the standard case of Rota calculus realisation. Nevertheless
the overall picture and system of statements depending only on GHW algebra
is the same modulo some automatic replacements in formulas demonstrated in
the sequel. The large part of that kind of job was already done in \cite{2,3}.

The aim of this presentation is to give a general picture of the
algebra of linear operators on polynomial algebra. The picture
that emerges discloses the fact that any $\psi $-representation of
finite operator calculus or equivalently - any $\psi
$-representation of GHW algebra makes up an example of the
algebraization of the analysis with generalized differential
operators \cite{11} acting on the algebra of polynomials.

We shall delimit all our considerations to the algebra $P$ of polynomials or sometimes to
the algebra of formal series. Therefore the distinction between
difference and differentiation operators disappears. All linear
operators on $P$ are both difference and differentiation operators
if the degree of differentiation or difference operator is
unlimited.

If all this is extended to Markowsky $Q$-umbral calculus \cite{11}
then many of the results of $\psi $-calculus may be extended to
$Q$-umbral calculus \cite{11}. This is achieved under the almost
automatic replacement of $\{D,\, \hat {x}, id\}$ generators of GHW
or their $\psi $-representation $\{\partial _{\psi},\hat
{x}_{\psi},id\}$ by their $Q$-representation correspondents
$\{Q,\hat {x}_{Q} ,id\}$ - see definition 2.5.

\vspace{0.1cm}
\noindent The article is supplemented by the short indicatory glossaries of
terms and notation used by Ward \cite{1}, Viskov
\cite{6},\cite{7}, Markowsky\cite{11}, Roman\cite{27}-\cite{30} on
one side and the Rota-oriented\cite{8}\-cite{10} notation on the
other side \cite{2}\,cite{3}.

\section{Primary definitions, notation and general observations}

In the following we shall consider the algebra $P$ of polynomials
$P = ${\bf F}[x] over the field {\bf F}  of characteristic zero. All
operators or functionals studied here are to be understood as {\it linear}
operators on $P$. It shall be easy to see that they are always well defined.

Throughout the note while saying ``polynomial sequence $\left\{
{p_{n}} \right\}_{0}^{\infty}  $'' we mean\\ deg $p_{n}= n$;\;$n
\geq 0$ and we adopt also the convention that deg $p_{n} < 0$  iff
$p_{n} \equiv 0$.

Consider $\Im $ - the family of functions` sequences (in conformity with
Viskov \cite{6},\cite{7},\cite{2} notation ) such that:\\
$\Im = \{\psi;R \supset \left[ {a,b} \right]\;;\;q \in \left[ {a,b}
\right]\;;\;\psi \left( {q} \right):Z \to F\;;\;\psi _{0} \left( {q} \right)
= 1\;;\;\psi _{n} \left( {q} \right) \ne 0;\;\psi _{ - n} \left( {q} \right)
= 0;\;n \in N\}$.\\
We shall call $\psi = \left\{ {\psi _{n} \left( {q} \right)} \right\}_{n
\ge 0} $ ; $\psi _{n} \left( {q} \right) \ne 0$; $n \ge 0$ and $\psi _{0}
\left( {q} \right) = 1$ an admissible sequence. Let now $n_{\psi}$
denotes \cite{2,3}  

$$ 
n_{\psi}  \equiv \psi _{n - 1} \left( {q} \right)\psi _{n}^{-1} \left( {q} \right),n \geq 0 .
$$
Then  (note that for admissible $\psi $,  $0_{\psi}= 0$)
$$
n_{\psi}! \equiv \psi _{n}^{-1} \left( {q} \right) \equiv
n_{\psi}
\left( {n - 1} \right)_{\psi}
\left( {n - 2} \right)_{\psi}
\left( {n - 3} \right)_{\psi}  ....
2_{\psi}
1_{\psi};
\quad
0_{\psi}!=1
$$

 $n_{\psi} ^{\underline {k}}  = n_{\psi}  \left( {n - 1} \right)_{\psi}
...\left( {n - k + 1} \right)_{\psi}  $,  $\left(
{{\begin{array}{*{20}c}
 {n} \hfill \\
 {k} \hfill \\
\end{array}} } \right)_{\psi}  \equiv \frac{{n_{\psi} ^{\underline {k}}
}}{{k_{\psi}  !}}  and  {\bf {\exp}}_{\psi}  \{ y\} =
\sum\limits_{k = 0}^{\infty}  {\frac{{y^{k}}}{{k_{\psi}  !}}} .$

\begin{defn}\label{deftwoone}
Let $\psi $ be admissible. Let $\partial _{\psi}  $ be the linear operator
lowering degree of polynomials by one defined according to $\partial _{\psi
} x^{n} = n_{\psi}  x^{n - 1}\;;\; n \geq 0$. Then $\partial _{\psi}  $
is called the $\psi $-derivative.
\end{defn}

\begin{rem}\label{remtwoone} {\em
a) For any rational function R the corresponding factorial $
{R\left( {q^{n}} \right)!}$  of the sequence $ R(q^{n})$ is
defined naturally  \cite{1,2} as it is defined for $ n_{\psi}$
sequence , i.e. : $R(q^{n})! = R(q^{n})R(q^{n-1})...R(q^{1})$

The choice $\psi _{n} \left( {q}
\right)$=$\left[ {R\left( {q^{n}} \right)!} \right]^{ - 1}$ and
$R\left( {x} \right) = \frac{{1 - x}}{{1 - q}}$ results in the
well known $q$-factorial\; $n_{q} ! = n_{q} \left( {n - 1}
\right)_{q} !;\quad 1_{q} ! = 0_{q} ! = 1$ while the $\psi
$-derivative $\partial _{\psi}  $ becomes now ($n_{\psi} = n_{q}
$) the Jackson's derivative \cite{24,25,26,2,3} $\partial _{q} $:
\begin{center}
$\left( {\partial _{q} \varphi}  \right)\left( {x} \right) = \frac{{\varphi
\left( {x} \right) - \varphi \left( {qx} \right)}}{{\left( {1-q} \right)x}}$.
\end{center}
b) Note also that if $\psi = \left\{ {\psi _{n} \left( {q}
\right)} \right\}_{n \geq 0} $ and $\varphi = \left\{ {\varphi
_{n} \left( {q} \right)} \right\}_{n \ge 0} $ are two admissible
sequences then [$\partial _{\psi}  $ , $\partial _{\varphi}$]$ =
0$ iff $\psi = \varphi$.Here [,] denotes the commutator of
operators.}
\end{rem}

\begin{defn}\label{deftwotwo}
Let $E^{y}\left( {\partial _{\psi} }  \right) \equiv exp_{\psi} \{
y\partial _{\psi}  \} = \sum\limits_{k = 0}^{\infty}
{\frac{{y^{k}\partial _{\psi}  ^{k}}}{{k_{\psi}!}}} $.\;
$E^{y}\left( {\partial _{\psi} }\right)$ is called the generalized
translation operator.
\end{defn}

\begin{note}\label{notetwoone}
{\em \cite{2,3}\\
$E^{a}\left( {\partial _{\psi} }  \right) f(x) \equiv f(x+_{\psi}  a)\;;\;
(x +_{\psi} a)^n \equiv E^{a}\left( {\partial _{\psi} }
\right) x^{n}\;;\;E^{a}\left( {\partial _{\psi} }  \right) f
= \sum \limits_{n \geq 0} {\frac{{a^{n}}}{{n_{\psi}  !}}} \partial _{\psi
}^{n}f$;\\
and in general $(x +_{\psi} a)^n \ne (x +_{\psi} a)^{n-1} (x +_{\psi} a)$.\\
Note also \cite{1} that in general $(1 +_{\psi} (-1))^{2n+1} \ne
0$ ; $n \geq 0$ though $(1 +_{\psi} (-1))^{2n} = 0$ ; $n \geq 1$.}
\end{note}

\begin{note}\label{notetwotwo}
{\em \cite{1}\\
$\exp_{\psi}  \left( {x + _{\psi}  y} \right) \equiv
exp_{\psi}  \{ x\} \exp_{\psi}  \{ y\}$ - while in general $\exp_{\psi}
\{ x + y\} \ne \exp_{\psi}  \{ x\} \exp_{\psi}  \{ y\} $.}
\end{note}
Possible consequent utilization of the identity $\exp_{\psi}  \left( {x +
_{\psi}  y} \right) \equiv \exp_{\psi}  \{ x\} \exp_{\psi}  \{ y\} $ is
quite encouraging. It leads among others to ``$\psi $-trigonometry''
either $\psi $-elliptic or $\psi $-hyperbolic via introducing $\cos_{\psi},\,
\sin_{\psi}$ \cite{1}, $\cosh_{\psi}\,,\,\sinh_{\psi}  $ or in general $\psi
$-hyperbolic functions of m-th order $\left\{ {h_{j}^{\left( {\psi}
\right)} \left( {\alpha}  \right)} \right\}_{j \in Z_{m}}  $defined
according to \cite{12}
\[
R \ni \alpha \to \;h_{j}^{(\psi)} \left( {\alpha}  \right) =
\frac{{1}}{{m}}\sum\limits_{k \in Z_{m}}  {\omega ^{ -
kj}\exp_{\psi}\left\{ {\omega ^{k}\alpha}  \right\}\,;\;j \in
Z_{m}},\;\omega = \exp\left\{ {i\frac{{2\pi} }{{m}}} \right\}.
\]
where $1 < m \in N$ and $Z_{m} = \{ 0,1,...,m-1\}$.

\begin{defn}\label{deftwothree}
A polynomial sequence $\left\{ {p_{n}}  \right\}_{o}^{\infty}  $ is of
$\psi${\it -binomial} type if it satisfies the recurrence
\[
E^{y}\left( {\partial _{\psi} }  \right) p_{n} \left( {x} \right) \equiv
p_{n} \left( {x +_{\psi}  y} \right) \equiv \sum\limits_{k \geq 0} {\left(
{{\begin{array}{*{20}c}
 {n} \hfill \\
 {k} \hfill \\
\end{array}} } \right)} _{\psi}  p_{k} \left( {x} \right)p_{n - k} \left(
{y} \right).
\]
\end{defn}

\noindent Polynomial sequences of $\psi $-binomial type \cite{2,3} are known
to correspond in one-to-one manner to special generalized
differential operators $Q$, namely to those $Q = Q\left( {\partial
_{\psi} } \right)$ which are $\partial _{\psi}  $-shift invariant
operators \cite{2,3}. We shall deal in this note mostly with this
special case,i.e. with $\psi $-umbral calculus. However before to
proceed let us supply a basic information referring to this
general case of $Q$-umbral calculus.

\begin{defn}\label{deftwofour}
Let $P = {\bf F}$[x]. Let $Q$ be a linear map $Q\;:\;P \to P$ such that:\\
$\forall_{p \in P}$ deg $(Qp) = ($deg $p) -1$ (with the convention deg $p =
-1$ means $p = const = 0$). $Q$ is then called a generalized
difference-tial operator \cite{11} or Gel`fond-Leontiev \cite{7} operator.
\end{defn}

\vspace{0.1cm}
\noindent From the above definitions we infer that the following holds.

\begin{obs}\label{obstwoone} {\em
Let $Q$ be as in Definition~\ref{deftwofour}. Let $Qx^{n}=
\sum\limits_{k = 1}^{n} {b_{n,k}}  x^{n - k}$ where $b_{n,1} \ne
0$ of course. Without loose of generality take $b_{1,1} = 1$. Then
$\exists\;{\left\{ {q_{k}}  \right\}_{k \geq 2}  \subset {\bf F}}$
and there exists  admissible $\psi $ such that
\begin{equation} \label{eq1}
Q = \partial _{\psi} + \sum\limits_{k \geq 2} {q_{k} \partial
_{\psi} ^{k}}
\end{equation}
if and only if
\begin{equation} \label{eq2}
b_{n,k} = \left( {{\begin{array}{*{20}c}
 {n} \hfill \\
 {k} \hfill \\
\end{array}} } \right)_{\psi}
b_{k,k}; \quad  n \geq k \geq 1;\;b_{n,1} \ne 0;\;b_{1,1} = 1.
\end{equation} 
If $\left\{ {q_{k}} \right\}_{k \geq 2} $ and an admissible $\psi$ exists then these are unique.} 
\end{obs}

\begin{n} {\em
In  case (\ref{eq2}) is true we shall write : $Q = Q\left(
{\partial _{\psi} }  \right)$ because then and only then the
generalized differential operator  Q  is a series in powers of
$\partial _{\psi}  $ .}
\end{n}

\begin{rem} {\em
Note that operators of the (\ref{eq1}) form constitute a group
under superposition of formal power series (compare with the
formula (S) in \cite{13}). Of course not all generalized
difference-tial operators satisfy (\ref{eq1}) i.e. are series just
only in corresponding $\psi $-derivative $\partial _{\psi}  $ (see
Proposition~\ref{propthreeone} ). For example \cite{14} let
$Q=\frac{1}{2}D \hat {x} D- \frac{1}{3}D^{3}$.\; Then $Qx^{n} =
\frac{1}{2}n^{2}x^{n - 1} - \frac{1}{3}n^{\underline {3}} x^{n -
3}$ so according to Observation~\ref{obstwoone} $n_{\psi} =
\frac{1}{2}n^{2}$ and there exists no admissible $\psi $ such that
$Q = Q\left( {\partial _{\psi} } \right)$. Here  $ \hat {x}$
denotes the operator of multipliation by  x  while $n^{\underline
{k}}$ is a special case of  $n_{\psi} ^{\underline {k}}$ for the
choice  $n_{\psi}$ = n.}
\end{rem}

\begin{obs} \label{obstwotwo} {\em
From theorem 3.1 in \cite{11} we infer that generalized differential operators
give rise to subalgebras $\sum_{Q}$ of linear maps (plus zero
map of course) commuting with a given generalized difference-tial operator
$Q$. The intersection of two different algebras $\sum_{Q_{1}}$ and
$\sum_{Q_{2}}$ is just zero map added.}
\end{obs}
\vspace{0.1cm} 
\noindent The importance of the above Observation~\ref{obstwotwo} as well as the
definition below may be further fully appreciated in the context of the
Theorem~\ref{thmtwoone} and the Proposition~\ref{propthreeone} to come.

\begin{defn} \label{deftwofive}
Let $\left\{ {p_{n}}  \right\}_{n \geq 0} $ be the normal
polynomial sequence \cite{11} ,i.e. $p_{0} \left( {x} \right)= 1$
and $p_{n} \left( {0} \right) = 0\;;\;n\geq 1$. Then we call it
the $\psi $-basic sequence of the generalized difference-tial
operator $Q$ if in addition $Q\,p_{n} = n_{\psi}p_{n - 1}$.
Parallely we define a linear map $\hat {x}_{Q} $: $P \to P$ such
that $\hat {x}_{Q} p_{n} =\frac{{\left( {n + 1} \right)}}{{\left(
{n + 1} \right)_{\psi} } }p_{n + 1} ;\quad n \geq 0$. We call the
operator $\hat {x}_{Q} $ the {\it dual} to $Q$ operator.
\end{defn}
When $Q = Q\left( {\partial _{\psi} }  \right) = \partial
_{\psi}  $ we write for short: $\hat {x}_{Q\left( {\partial _{\psi}}  \right)}
\equiv \hat {x}_{\partial _{\psi}}  \equiv \hat {x}_{\psi}$ (see: Definition
\ref{deftwonine}).\\
Of course [$Q,\hat {x}_{Q} $]$ = id$ therefore $\{Q,{\hat {x}}_{Q},
id\}$ provide us with a continuous family of generators of GHW in -
as we call it - $Q$-representation of Graves-Heisenberg-Weyl algebra.\\
In the following we shall restrict to special case of generalized
differential operators $Q$, namely to those $Q = Q\left( {\partial
_{\psi} }  \right)$ which are $\partial _{\psi}  $-shift invariant operators
\cite{2,3} (see: Definition \ref{deftwosix}).

\vspace{0.1cm}
\noindent  Let us start with appropriate $\psi $-Leibniz rules for
corresponding $\psi $-derivatives.\\
$\psi${\bf -Leibniz rules:}

\vspace{0.2cm}

\noindent It is easy to see that the following hold for any formal series $f$ and
$g$:\\
for $\partial _{q} $:\;\; $\partial _{q} \left( {f \cdot g} \right) =
\left({\partial _{q} f} \right) \cdot g + \left( {\hat {Q}f} \right) \cdot
\left({\partial _{q} g} \right)$, where $\left( {\hat {Q}f} \right) \left( {x}
\right) = f \left( {qx} \right)$;\\
for $\partial _{R} = R\left( {q\hat {Q}} \right)\partial _{0}
$:\;\; $\partial _{R} (f \cdot g) (z) = R\left( {q\hat {Q}}
\right)\{
(\partial _{0} f) (z)  \cdot  g(z) + f(0) (\partial _{0} g) (z)\}$\\
where - note - $R\left( {q\hat {Q}} \right)x^{n - 1}= n_{R} x^{n - 1}$
; $(n_{\psi}  = n_{R} =n_{R(q)} = R\left( {q^{n}} \right))$ and finally\\
for $\partial _{\psi} =\hat {n}_{\psi}  \partial _{0} $:
\[
\partial _{\psi} (f \cdot g) (z) =\hat {n}_{\psi} \{ (\partial _{o}
f) (z)  \cdot  g(z) + f(0) (\partial _{0} g) (z)\}
\]
\vspace{0.1cm}
\noindent where $\hat {n}_{\psi}  x^{n - 1}= n_{\psi}  x^{n - 1}$ ; $n \geq
1$.

\begin{ex} {\em
Let $Q\left( {\partial _{\psi} }  \right)= D \hat {x} D$, where
$\hat {x} f(x)= x f(x)$ and $D = \frac{d}{dx}$.
Then $\psi = \left\{ {\left[ {\left( {n^{2}} \right)!} \right]^{ - 1}}
\right\}_{n \geq 0} $ and $Q = \partial _{\psi}  $.
Let $Q\left( {\partial _{\psi} }  \right) R(q \hat {Q}) \,
\partial _{0} \equiv \partial _{R} $. Then $\psi = \left\{
{\left[ {R\left( {q^{n}} \right)!} \right]^{ - 1}} \right\}_{n \geq 0} $ and
$Q = \partial _{\psi} \equiv \partial_{R} $. Here $R(z)$ is
any formal Laurent series; $\hat {Q}f(x) = f(qx)$
and $n_{\psi} =R(q^{n})$. $\partial _{0} $ is $q = 0$
Jackson derivative which as a matter of fact - being a difference operator
is the differential operator of infinite order at the same time:
\begin{equation}
\label{eq3}
\partial _{0} =
\sum\limits_{n = 1}^{\infty}  {\left( { - 1} \right)^{n + 1}\frac{{x^{n -
1}}}{{n!}}\frac{{d^{n}}}{{dx^{n}}}}.
\end{equation}

\vspace{0.1cm}
\noindent Naturally with the choice $\psi _{n} \left( {q} \right)=\left[
{R\left( {q^{n}} \right)!} \right]^{ - 1}$ and $R\left( {x}
\right) = \frac{{1 - x}}{{1 - q}}$ the $\psi $-derivative
$\partial _{\psi}  $ becomes the Jackson's derivative
\cite{24,25,26,2,3} $\partial _{q} $:
\[
\left( {\partial _{q} \varphi}  \right)\left( {x} \right) = \frac{{1-q\hat
{Q}}}{{\left( {1 - q} \right)}}\partial _{0} \varphi \left( {x} \right).
\]
}
\end{ex}

\vspace{0.1cm}
\noindent  The equivalent to (\ref{eq3}) form of Bernoulli-Taylor expansion one may
find \cite{15} in {\it Acta Eruditorum} from November 1694 under the name
``{\it series univeralissima''}.

\vspace{0.1cm}
\noindent (Taylor`s expansion was presented in his ``Methodus incrementorum directa
et inversa'' in 1715 - edited in London).

\begin{defn} \label{deftwosix}
Let us denote by $End(P)$ the algebra of all linear operators acting on the
algebra $P$ of polynomials. Let
$$
\sum\nolimits_{\psi} = \{ T \in End (P);\;\forall \;\alpha \in \bf
F;\; \left[ {T,E^{\alpha} \left( {\partial _{\psi} } \right)}
\right] = 0 \}.
$$
\vspace{0.1cm}
\noindent Then $\sum_{\psi}$ is a commutative subalgebra of $End (P)$ of
${\bf F}$-linear operators. We shall call these operators
$T:\;\partial _{\psi}  $-shift invariant operators.
\end{defn}
\vspace{0.1cm}
\noindent We are now in a position to define further basic objects of ``$\psi $-umbral
calculus'' \cite{2,3}.

\begin{defn} \label{deftwoseven}
Let $Q\left( {\partial _{\psi} }  \right):P \to P;$ the linear operator
$Q\left( {\partial _{\psi} }  \right)$ is a $\partial _{\psi}$-{\it delta}
operator iff
\begin{enumerate}
\renewcommand{\labelenumi}{\em \alph{enumi})}
\item $Q\left( {\partial _{\psi} }  \right)$ is $\partial _{\psi}  $ - shift
invariant;
\item $Q\left( {\partial _{\psi} }  \right)\left( {id} \right) =
const \ne 0$  where  id(x)=x.
\end{enumerate}
\end{defn}

\vspace{0.1cm}
\noindent The strictly related notion is that of the $\partial _{\psi}  $-basic
polynomial sequence:

\begin{defn} \label{deftwoeight}
Let $Q\left( {\partial _{\psi} }  \right):P \to P;$ be the $\partial _{\psi
} $-\textit{delta} operator. A polynomial sequence
$\left\{ {p_{n}}  \right\}_{n \ge 0} $;
deg \textit{p}$_{n}$\textit{= n} such that:
\begin{enumerate}
\renewcommand{\labelenumi}{\em \arabic{enumi})}
\item $p_{0} \left( {x} \right) = 1$;
\item $p_{n} \left( {0} \right) = 0$; $n > 0$;
\item $Q\left( {\partial
_{\psi} }  \right)p_{n} = n_{\psi}  p_{n - 1} $ ,$\partial _{\psi}
$-delta operator$ Q\left( {\partial_ \psi} \right)$is called the
$\partial _{\psi}  $-basic polynomial sequence of the $\partial
_{\psi}  $-delta operator.
\end{enumerate}
\end{defn}

\begin{ident} {\em
It is easy to see that the following identification takes place: $\partial
_{\psi}  $-delta operator $Q\left( {\partial _{\psi} }  \right) = \partial
_{\psi}  $-shift invariant generalized differential operator $Q$.
Of course not every generalized differential operator might be considered
to be such.}
\end{ident}

\vspace{0.2cm}
\noindent {\bf Note 2.3.}: Let $\Phi \left( {x;\lambda}  \right) =
\sum\limits_{n \geq 0} {\frac{{\lambda ^{n}}}{{n_{\psi}  !}}p_{n}
\left( {x} \right)} $ denotes the $\psi $-exponential generating
function of the $\partial _{\psi}$-basic polynomial sequence
$\left\{ {p_{n}} \right\}_{n \geq 0} $ of the $\partial_{\psi}
$-delta operator $Q \equiv Q\left( {\partial _{\psi}} \right)$ and
let $\Phi \left( {0;\lambda}  \right) = 1$. Then $Q \Phi \left(
{x;\lambda} \right) = \lambda \Phi \left( {x;\lambda} \right)$ and
$\Phi $ is the unique solution of this eigenvalue problem. If in
addition (2.2) is satisfied  then there exists such an admissible
sequence $\varphi $ that $\Phi \left( {x;\lambda} \right)=
\exp_{\varphi} \left\{\lambda x
 \right\} $  (see Example 3.1) .

\vspace{0.2cm}
\noindent The notation and naming established by Definitions
\ref{deftwoseven} and \ref{deftwoeight} serve the target to
preserve and to broaden simplicity of Rota`s finite operator
calculus also in its extended ``$\psi $-umbral calculus'' case
\cite{2,3}. As a matter of illustration of such notation
efficiency let us quote after \cite{2} the important
Theorem~\ref{thmtwoone} which might be proved using the fact that
$\forall \; Q\left( {\partial _{\psi} }  \right) \quad \exists !$
invertible $S \in \Sigma_{\psi}$ such that $Q\left( {\partial
_{\psi} } \right)= \partial _{\psi} S$. ( For
Theorem~\ref{thmtwoone} see also Theorem 4.3. in \cite{11}, which
holds for operators, introduced by the
Definition~\ref{deftwofive}). Let us define at first what follows.

\begin{defn} \label{deftwonine} {\em (compare with (17) in \cite {7})}\\
The Pincherle $\psi $-derivative is the linear map ' : $\Sigma
_{\psi} \quad \to  \quad \Sigma _{\psi} $;
\begin{center}
 $T\;$ ' = $T\;\hat {x}_{\psi}  $ \textbf{-} $\hat {x}_{\psi}  T \equiv
 $\textbf{[}$T$, $\hat {x}_{\psi} $\textbf{]}
\end{center}
\vspace{0.1cm}
\noindent where the linear map $\hat {x}_{\psi}  :P \to P;$ is defined in
the basis $\left\{ {x^{n}} \right\}_{n \ge 0} $ as follows
\[
\hat {x}_{\psi}  x^{n} = \frac{{\psi _{n + 1} \left( {q}
\right)\left( {n + 1} \right)}}{{\psi _{n} \left( {q}
\right)}}x^{n + 1} = \frac{{\left( {n + 1} \right)}}{{\left( {n +
1} \right)_{\psi} } }x^{n + 1};\quad n \ge 0
\]
\end{defn}

\vspace{0.2cm}
\noindent The following theorem is true \cite{2}

\begin{thm}\label{thmtwoone}
{\em ($\psi $-Lagrange and $\psi$-Rodrigues formulas \cite{33,10,11,22,2})}\\
Let $\left\{ {p_{n} \left( {x} \right)} \right\}_{n = 0}^{\infty}  $ be
$\partial _{\psi}  $-basic polynomial sequence of the $\partial _{\psi}
$-delta operator $Q\left( {\partial _{\psi} }  \right)$.\\
Let $Q\left( {\partial _{\psi} }  \right) = \partial _{\psi}
S.$ Then for $n>0$:

\begin{enumerate}
\renewcommand{\labelenumi}{\em (\arabic{enumi})}
\item\label{one} $p_{n}(x) = Q\left( {\partial _{\psi} }
\right)$\textbf{'} $S^{-n-1}\;${x}$^{n}$ ;

\item\label{two} $p_{n}(x) = S^{-n}${x}$
^{n} - \frac{{n_{\psi} } }{{n}}$ ($S^{ -
n}\;$)\textbf{'}{x}$^{n-1};$

\item\label{three} $p_{n}(x) = \frac{{n_{\psi} } }{{n}}\hat {x}_{\psi}
S^{ - n}${x}$^{n-1}$;

\item\label{four} $p_{n}(x) = \frac{{n_{\psi} } }{{n}}\hat {x}_{\psi}
(Q\left( {\partial _{\psi} }  \right)$\textbf{'} )$^{-1}
p_{n-1}(x)$  {\em ($\leftarrow$ Rodrigues $\psi $-formula )}.
\end{enumerate}
\end{thm}

\vspace{0.2cm}
\noindent For the proof one uses typical properties of the Pincherle $\psi
$-derivative \cite{2}.

\begin{obs} \label{obstwothree} {\em [2,3]\\
The triples $\{\partial _{\psi},\hat {x}_{\psi}, id \}$ for any
admissible $\psi $-constitute the set of generators of the $\psi
$-labelled representations of Graves-Heisenberg-Weyl (GHW) algebra
\cite{16,17,18}. Namely, as easily seen [$\partial _{\psi},\hat {x}_{\psi}  $]
$= id$. (compare with Definition~\ref{deftwofive})}
\end{obs}

\begin{obs} {\em
In view of the Observation~\ref{obstwothree} the general Leibnitz rule in
$\psi $-representation of Graves-Heisenberg-Weyl algebra may be written
(compare with 2.2.2 Proposition in \cite{17}) as follows
\begin{equation}
\label{eq4}
\partial _{\psi} ^{n}
\quad
\hat {x}_{\psi} ^{m} \quad =
\quad
\sum\limits_{k \ge 0} {\left( {{\begin{array}{*{20}c}
 {n} \hfill \\
 {k} \hfill \\
\end{array}} } \right)\left( {{\begin{array}{*{20}c}
 {m} \hfill \\
 {k} \hfill \\
\end{array}} } \right)k!}\,
\hat {x}_{\psi} ^{m - k}\,
\partial _{\psi} ^{n - k}.
\end{equation}
\vspace{0.2cm}
\noindent One derives the above $\psi $-Leibnitz rule from $\psi$-Heisenberg-Weyl
exponential commutation rules exactly the same way as in
$\{D,\hat {x},id \}$ GHW representation - (compare with 2.2.1
Proposition in \cite{17} ). $\psi $-Heisenberg-Weyl exponential commutation
relations read:
\begin{equation} \label{eq5}
{\bf { exp}} \{t \partial _{\psi} \} {\bf {exp}} \{a
\hat{x}_{\psi} \} = {\bf { exp}} \{at \} {\bf {\ exp}} \{a
\hat{x}_{\psi}  \} {\bf {\it exp}} \{t \partial _{\psi} \}.
\end{equation}  }
\end{obs}
\vspace{0.2cm}
\noindent To this end let us introduce a pertinent $\psi $-multiplication
$ \ast _{\psi} $ of functions as specified below.

\begin{n} \label{notetwothree} {\textrm{  }\\
x $\ast _{\psi}$ x$^{n} = \hat {x}_{\psi}($x$^{n}) =
\frac{{\left( {n + 1} \right)}}{{\left( {n + 1} \right)_{\psi} } }x^{n +
1};\quad n \geq 0$ \; hence x $\ast _{\psi} 1 = 1_{\psi} $ x $\not
\equiv$ x \\
x$^{n} \ast _{\psi}$ x $= \hat {x}_{\psi} ^{n} ( $ x $ ) =
\frac{{\left( {n +1} \right)!}}{{\left( {n + 1} \right)_{\psi}  !}}x^{n + 1}
;\quad n \geq 0$
hence $1 \ast _{\psi}$ x $= 1_{\psi}$ x $\not \equiv $ x therefore\\
x $ \ast _{\psi} \alpha 1 = \alpha 1 \ast _{\psi} $ x $=$ x $ \ast
_{\psi} \alpha = \alpha \ast _{\psi} $ x $= \alpha  1_{\psi}$ x
and $\forall$ x,$\alpha \in {\bf F};\; f( x) \ast _{\psi} $ x$^{n}
= f(\hat {x}_{\psi})$x$^{n}$. }
\end{n}

\vspace{0.2cm}
\noindent For $k \ne n $ \; x$^{n} \ast _{\psi} $ x$^{k} \ne$ x$^{k} \ast _{\psi} $
x$^{n}$ as well as x$^{n} \ast _{\psi} $ x$^{k} \ne$ x$^{n+ k}$ - in
general i.e. for arbitrary admissible $\psi $; compare this
with $($x $+_{\psi} a)^n \ne ($ x $+_{\psi} a)^{n-1}( $x $+_{\psi} a)$.\\
In order to facilitate in the future formulation of observations accounted
for on the basis of $\psi $-calculus representation of GHW algebra we shall
use what follows.

\begin{defn}
With Notation \ref{notetwotwo} adopted let us define the $ \ast _{\psi} $
powers of x according to

 $x^{n \ast _{\psi} } \equiv $ x $ \ast _{\psi} x^{\left( {n - 1} \right)
\ast _{\psi} } = \hat {x}_{\psi} (x^{\left( {n - 1} \right)\ast _{\psi
}} ) = $ x $ \ast _{\psi} $ x $ \ast _{\psi} $ ... $ \ast _{\psi} $ x
$=\frac{n!}{n_{\psi}  !}x^{n};\quad n \geq 0$.
\end{defn}
\vspace{0.2cm}
\noindent Note that $x^{n\ast _{\psi} }  \ast _{\psi} x^{k\ast _{\psi}
} = \frac{{n!}}{{n_{\psi}  !}} x^{\left( {n + k} \right)\ast _{\psi} }
\ne x^{k\ast _{\psi} }  \ast _{\psi} x^{n\ast _{\psi} }  =
\frac{{k!}}{{k_{\psi}  !}}x^{\left( {n + k} \right)\ast _{\psi} } $ for
$k  \ne n$ and $x^{0\ast _{\psi} }=1$.\\
\vspace{0.2cm}
\noindent This noncommutative $\psi $-product $ \ast _{\psi}$ is deviced so as to
ensure the following observations.

\begin{obs} \label{obstwofive} {\em \textrm{                      }
\begin{enumerate}
\renewcommand{\labelenumi}{\em \alph{enumi})}
\item $\partial _{\psi}  x^{n\ast _{\psi} } = n x^{\left( {n - 1}
 \right)\ast _{\psi} } $;\; $n \ge 0$
\item  ${\bf \exp}_{\psi}  [\alpha {\it x}]  \equiv {\bf\exp}
\{\alpha \hat {x}_{\psi}  $\}{\bf 1}
\item ${\bf\exp} [\alpha x]  \ast _{\psi} 
({\bf\exp}_{\psi} \{\beta \hat {x}_{\psi}  \}{\bf 1}) = ({\bf
\exp}_{\psi}  \{[\alpha +\beta ]\hat {x}_{\psi}  $\}){\bf 1}
\item $\partial _{\psi} (x^{k} \ast _{\psi} \quad x^{n\ast _{\psi} } )
= (D x^{k}) \ast _{\psi} x^{n \ast _{\psi} } + x^{k} \ast _{\psi}
(\partial _{\psi}  x^{n\ast _{\psi} })$ hence
\item \label{e}
$\partial _{\psi} ( f \ast _{\psi} g) = ( Df) \ast _{\psi} g  + f
\ast _{\psi} (\partial _{\psi} g)$ ; $f,g$  - formal series
\item \label{f}
$f( \hat {x}_{\psi}) g (\hat {x}_{\psi} )$ {\bf 1} $= f(x) \ast
_{\psi} \tilde {g} (x)$ ; $\tilde {g} (x) = g(\hat {x}_{\psi})${\bf 1}.
\end{enumerate} }
\end{obs}
\vspace{0.2cm}
\noindent Now the consequences of Leibnitz rule (e) for difference-ization of the
product are easily feasible. For example the Poisson $\psi $-process
distribution ${\pi}_{m}(x) = \frac{1}{N(\lambda , x)}$p$_m(x);\;
\sum\limits_{m \geq 0}\,$p$_{m}(x) = 1$ is determined by
\begin{equation} \label{eq6}
\textrm{p}_{m}(x) = \frac{{\left( {\lambda x} \right)^{m}}}{{m!}} \ast
_{\psi} {\bf exp}_{\psi} [-\lambda x]
\end{equation}
\vspace{0.1cm}
\noindent which is the unique solution (up to a constant factor) of the
 $\partial _{\psi}
$-difference equations systems
\begin{equation} \label{eq7}
\partial _{\psi} \textrm{p}_{m}(x) + \lambda  \textrm{p}_{m}(x) = \lambda
\textrm{p}_{m-1}(x) \; m > 0\;;\; \partial _{\psi} \textrm{p}_{0} (x) =
- \lambda  \textrm{p}_{0} (x)
\end{equation}
\vspace{0.2cm}
\noindent Naturally $N(\lambda , x) = {\bf\exp} [\lambda x] \ast
_{\psi} {\bf\exp}_{\psi} [-\lambda x$].

\vspace{0.2cm}
\noindent As announced - the rules of $\psi $ -product $ \ast _{\psi} $ are
accounted for - as a matter of fact - on the basis of $\psi
$-calculus representation of GHW algebra. Indeed,it is enough to
consult Observation~\ref{obstwofive} and to introduce $\psi
$-Pincherle derivation $\hat {\partial} _{\psi}  $ of series in
powers of the symbol $\hat {x}_{\psi}  $ as below. Then the
correspondence between generic relative formulas turns out
evident.

\begin{obs} {\em
Let $\hat {\partial} _{\psi}  \equiv \frac{{\partial} }{{\partial \hat
{x}_{\psi} } }$ be defined according to $\hat {\partial} _{\psi}
f(\hat {x}_{\psi}) = [\partial _{\psi},f(\hat {x}_{\psi})]$. Then
$\hat {\partial} _{\psi} \hat {x}_{\psi} ^{n} = n \hat {x}_{\psi} ^{n - 1}\;
;\; n \geq 0$ and $\hat {\partial} _{\psi}  \hat {x}_{\psi} ^{n}$ {\bf 1}
$=\partial_{\psi}  x^{n\ast _{\psi} }$ hence {\bf [}$\hat {\partial}_{\psi}
f(\hat {x}_{\psi})${\bf ]1 } $=\partial _{\psi}f(x)$ where $f$ is a
formal series in powers of $\hat {x}_{\psi}  $ or equivalently in $ \ast
_{\psi}$ powers of $x$.}
\end{obs}

\vspace{0.2cm}
\noindent As an example of application note how the solution of \ref{eq7} is obtained
from the obvious solution ${\mathbf p}_{m}(\hat {x}_{\psi})$ of the $\hat
{\partial} _{\psi}  $-Pincherle differential equation \ref{eq8} formulated
within G-H-W algebra generated by $\{\partial _{\psi},\hat {x}_{\psi},id \}$
\begin{equation} \label{eq8}
 \hat {\partial} _{\psi} {\mathbf p}_{m} (\hat {x}_{\psi}) + \lambda
 {\mathbf p}_{m}(\hat {x}_{\psi}) = \lambda {\mathbf p}_{m-1}(\hat {x}_{\psi}
 )\, m >0\;;\;\partial _{\psi} {\mathbf p}_{0} (\hat {x}_{\psi}  ) = -
\lambda {\mathbf p}_{0} (\hat {x}_{\psi}. )
\end{equation}
\vspace{0.2cm}
\noindent Namely : due to Observation~\ref{obstwofive} (f) p$_{m}(x) = {\mathbf p}
_{m}(\hat {x}_{\psi})${\bf 1}, where
\begin{equation}
{\mathbf p}_{m}(\hat {x}_{\psi}) = \frac{{\left( {\lambda \hat {x}_{\psi}
} \right)^{m}}}{{m!}}{\bf exp}_{\psi} [-\lambda \hat {x}_{\psi}].
\end{equation}

\section{ The general picture}

\noindent The general picture from the title above relates to the general picture of
the algebra $End(P)$ of operators on $P$ as in the following we shall
consider the algebra $P$ of polynomials $P$ = {\bf F}[x]
over the field {\bf F}  of characteristic zero.

\vspace{0.2cm}
\noindent We shall draw an over view picture of the situation distinguished
by possibility to develop umbral calculus for {\it any} polynomial
sequences $\left\{ {p_{n}}  \right\}_{0}^{\infty}  $ instead of
those of traditional binomial type only.

\vspace{0.2cm}
\noindent In 1901 it was proved \cite{19} that every linear operator mapping $P$ into
$P$ may be represented as infinite series in operators $\hat {x}$ and
{\it D}. In 1986 the authors of \cite{20} supplied the explicit expression for
such series in most general case of polynomials in one variable ( for many
variables see: \cite{21} ). Thus according to Proposition 1 from \cite{20}
one has:

\begin{prop} \label{propthreeone} {\em
Let $Q$  be a linear operator that reduces by one the degree of each
polynomial.
Let $\{ q_{n} \left( {\hat {x}} \right)\} _{n \ge 0} $ be an arbitrary
sequence of polynomials in the operator $\hat {x}$. Then $\hat {T} =
\sum\limits_{n \ge 0} {q_{n} \left( {\hat {x}} \right)} Q^{n}$
defines a linear operator that maps polynomials into polynomials.
Conversely, if $\hat {T}$ is linear operator that maps polynomials into
polynomials then there exists a unique expansion of the form
\[
\hat {T} = \sum\limits_{n \ge 0} {q_{n} \left( {\hat {x}} \right)} Q^{n}.
\]
}
\end{prop}

\vspace{0.1cm}
\noindent It is also  an easy exercise to prove the Proposition 2 from
\cite{20}:

\begin{prop} \label{propthreetwo} {\em
{\it  Let} $Q$ {\it be a linear operator that reduces by one the
degree of each polynomial. Let} $\{ q_{n} \left( {\hat {x}}
\right)\} _{n \ge 0} ${\it be an arbitrary sequence of polynomials
in the operator} $\hat {x}$. {\it Let} {\it a linear operator that
maps polynomials into polynomials be given by }
\begin{center}
$\hat {T} = \sum\limits_{n \ge 0} {q_{n} \left( {\hat {x}} \right)} Q^{n}$.
\end{center}
{\it Let} $P\left( {x;\lambda}  \right) = \sum\limits_{n \ge 0}
{q_{n} \left( {x} \right)} \lambda ^{n}$ {\it denotes indicator
of} $\hat {T}${\it. Then there exists a unique formal series}
$\Phi \left( {x;\lambda } \right)$;\;$\Phi \left( {0;\lambda}
\right) = 1$  {\it such that:}
\begin{center}
 $Q \Phi \left( {x;\lambda}  \right) = \lambda
\Phi \left( {x;\lambda}  \right).$
\end{center}
{\it Then also} $P\left( {x;\lambda} \right) = \Phi \left( {x;\lambda}
\right)^{ - 1}\hat {T}\Phi \left({x;\lambda}  \right)$. }
\end{prop}

\begin{ex} {\em
\vspace{0.2cm}
\noindent Note that \; $\partial _{\psi}  {\bf\exp}_{\psi} \{\lambda x \}=
\lambda \bf{\exp}_{\psi} \{\lambda x\};\; \bf{exp}_{\psi}
\left[ {x}
\right]\left| {_{x = 0}}  \right. = 1$. \;\;(*)\\
\vspace{0.1cm}
\noindent Hence for indicator of $\hat {T};\; \hat {T} = \sum\limits_{n \ge 0}
{q_{n} \left( {\hat {x}} \right)} \partial _{\psi} ^{n} $ we have:
\begin{center}
$P\left( {x;\lambda}  \right)= [{\bf{\exp}}_{\psi} \{\lambda {\it
x}\}]^{-1}\hat{T}\bf{\exp}_{\psi}  \{\lambda {\it
x}\}.\;\; $ (**)  
\end{center}
\vspace{0.1cm}
\noindent After choosing $\psi _{n} \left( {q} \right)=\left[ {n_{q} !}
\right]^{ - 1}$ we get $\bf {\exp}_{\psi}  \{{\it x}\} =\bf
{\exp}_{q} \{x\}$. In this connection note that $exp_{0} \left(
{x} \right) = \frac{{1}}{{1 - x}}$ and $exp(x)$ are mutual limit
deformations for $|x|<1$ due to:
\begin{center}
 $\frac{{exp_{0} \left( {z} \right) - 1}}{{z}} = exp_{o} \left( {z} \right)
\Rightarrow exp_{0} \left( {z} \right) = \frac{{1}}{{1 - z}} =
\sum\limits_{k = 0}^{\infty}  {z^{k}};\; |z|<1$ , i.e.
\end{center}
\[
exp\left( {x} \right)\mathrel{\mathop{\kern0pt\longleftarrow}\limits_{1
\leftarrow q}} exp_{q} \left( {x} \right) = \sum\limits_{n = 0}^{\infty}
{\frac{{x^{n}}}{{n_{q}
!}}\;\mathrel{\mathop{\kern0pt\longrightarrow}\limits_{q \to 0}}
\frac{{1}}{{1 - x}}}.
\]
\vspace{0.1cm}
\noindent Therefore corresponding specifications of (*) such as $exp_{0}
\left( {\lambda x} \right) = \frac{{1}}{{1 - \lambda x}}$ or
$\exp(\lambda x)$ lead to corresponding specifications of (**)
for divided difference operator $\partial _{0} $ and {\it D}
operator including special cases from \cite{20}.}
\end{ex}

\vspace{0.2cm}
\noindent To be complete let us still introduce \cite{2,3} an important operator $\hat
{x}_{Q\left( {\partial _{\psi} }  \right)} $ dual to $Q\left( {\partial
_{\psi} }  \right)$.

\begin{defn} \label{defthreeone} (see Definition~\ref{deftwofive})\\
Let $\left\{ {p_{n}}  \right\}_{n \ge 0} $ be the $\partial _{\psi}
$-basic polynomial sequence of the $\partial _{\psi}  $-delta
operator $Q\left( {\partial _{\psi} }  \right)$. A linear map $\hat
{x}_{Q\left( {\partial _{\psi} }  \right)} $: $P \to P$ ; $\hat {x}_{Q\left(
{\partial _{\psi} }  \right)}\, p_n =\frac{{\left( {n + 1} \right)}}{{\left( {n
+ 1} \right)_{\psi} } }p_{n + 1} ;\quad n \ge 0$ is called the operator
{\it dual} to $Q\left( {\partial _{\psi} }  \right)$.
\end{defn}

\begin{com} {\em
Dual in the above sense corresponds to adjoint in $\psi$-umbral calculus
language of linear functionals' umbral algebra (compare with Proposition
1.1.21 in \cite{22} ).}
\end{com}
\vspace{0.2cm}
\noindent It is now obvious that the following holds.
\begin{prop} {\it
Let $\{ q_{n} \left( {\hat {x}_{Q\left( {\partial _{\psi} }  \right)}}
\right)\} _{n \ge 0} $ be an arbitrary sequence of polynomials in the
operator $\hat {x}_{Q\left( {\partial _{\psi} }  \right)} $. Then $T =
\sum\limits_{n \ge 0} {q_{n} \left( {\hat {x}_{Q\left( {\partial _{\psi} }
\right)}}  \right)} Q\left( {\partial _{\psi} }  \right)^{n}$ defines a
linear operator that maps polynomials into polynomials. Conversely, if $T$
is linear operator that maps polynomials into polynomials then there exists
a unique expansion of the form
\begin{equation}
\label{eq9}
T = \sum\limits_{n \ge 0} {q_{n} \left( {\hat {x}_{Q\left( {\partial _{\psi
}}  \right)}}  \right)} Q\left( {\partial _{\psi} }  \right)^{n}.
\end{equation}   }
\end{prop}

\begin{com} {\em
The pair $Q\left( {\partial _{\psi} }  \right),\; \hat
{x}_{Q\left( {\partial _{\psi} }  \right)} $ of dual operators is expected
to play a role in the description of quantum-like processes apart from the
$q$-case now vastly exploited \cite{2,3}.}
\end{com}
\vspace{0.2cm}
\noindent Naturally the Proposition~\ref{propthreetwo} for $Q\left( {\partial _{\psi} }
\right)$ and $\hat {x}_{Q\left( {\partial _{\psi} }  \right)} $ dual
operators is also valid.\\
{\bf Summing up}: we have the following picture for $End(P)$ -
the algebra of all linear operators acting on the algebra $P$ of
polynomials.
\begin{center}
$Q(P) \equiv \bigcup\limits_{Q} \sum\nolimits_{Q} \subset End(P)$
\end{center}

\noindent and of course $Q(P) \ne End(P)$ where the subfamily $Q(P)$ (with zero map)
breaks up into sum of subalgebras $\sum\nolimits_{Q}$ according to
commutativity of these generalized difference-tial operators $Q$ (see
Definition~\ref{deftwofour} and Observation~\ref{obstwotwo}).
Also to each subalgebra $\sum\nolimits_{\psi}$ i.e. to each $Q\left(
{\partial _{\psi} }  \right)$ operator there corresponds its dual operator
$\hat {x}_{Q\left( {\partial _{\psi} }\right)} $
\[
\hat {x}_{Q\left( {\partial _{\psi} }  \right)} \notin
\sum\nolimits_{\psi}
\]

\noindent and both $Q\left( {\partial _{\psi} }  \right)$ \& $\hat {x}_{Q\left(
{\partial _{\psi} }  \right)} $ operators are sufficient to build up the
whole algebra $End(P)$ according to unique representation given by
(\ref{eq9}) including the $\partial
_{\psi}  $ and $\hat {x}_{\psi } $ case. Summarising: for any admissible
$\psi ${\it}  we have the following general statement.\\
{\bf General statement:}
\begin{center}
$End(P) = $[$\{\partial _{\psi}  $,$\hat {x}_{\psi}  $\}] =
[\{$Q\left( {\partial _{\psi} }  \right)$ , $\hat {x}_{Q\left( {\partial
_{\psi} }  \right)} $\}] = [\{$Q$ , $\hat {x}_{Q} $\}]
\end{center}
\vspace{0.1cm}
\noindent i.e. the algebra $End(P)$ is generated by any dual pair \{$Q$ ,
$\hat {x}_{Q} $\} including any dual pair \{$Q\left( {\partial _{\psi} }
\right)$ , $\hat {x}_{Q\left( {\partial _{\psi} }  \right)} $\} or
specifically by \{$\partial _{\psi}  $,$\hat {x}_{\psi}  $\} which in turn
is determined by a choice of any admissible sequence $\psi $.

\vspace{0.2cm}
\noindent As a matter of fact and in another words: we have bijective correspondences
between different commutation classes of $\partial _{\psi}  $-shift
invariant operators from $End(P)$, different abelian subalgebras
$\sum\nolimits_{\psi}$,\; distinct $\psi $-representations of GHW
algebra, different $\psi $-representations of the reduced incidence
algebra R(L(S)) - isomorphic to the algebra $\Phi _{\psi}  $ of $\psi
$-exponential formal power series \cite{2} and finally - distinct
$\psi $-umbral calculi \cite{7,11,14,23,2}. These bijective correspondences
may be naturally extended to encompass also $Q$-umbral calculi,
$Q$-representations of GHW algebra and abelian subalgebras
$\sum\nolimits_{Q}$.\\
(Recall: R(L(S)) is the reduced incidence algebra of L(S) where\\
L(S)=\{A; A$\subset $S; $|$A$|<\infty $\}; S is countable and (L(S);
$ \subseteq $) is partially ordered set ordered by inclusion \cite{10,2} ).

\vspace{0.2cm}
\noindent This is the way the Rota`s devise has been carried into effect. The devise
{\it ``much is the iteration of the few''} \cite{10} - much of the properties of
literally {\it all} polynomial sequences - as well as GHW algebra
representations - is the application of few basic principles of the
$\psi $-umbral difference operator calculus \cite{2}.\\
$\psi -$ {\bf Integration Remark :}\\
Recall : $\partial _{o} x^{n} = x^{n - 1}$. $\partial _{o} $ is identical
with divided difference operator. $\partial _{o} $ is identical with
$\partial _{\psi}  $ for $\psi = \left\{ {\psi \left( {q} \right)_{n}}
\right\}_{n \ge 0} $ ; $\psi \left( {q} \right)_{n} = 1$ ; $n \ge 0$ . Let
$\hat {Q}f(x)f(qx)$.\\
\vspace{0.2cm}
\noindent Recall also that there corresponds to the ``$\partial _{q} $
difference-ization'' the q-integration \cite{24,25,26} which is a right
inverse operation to ``$q$-difference-ization''. Namely
\begin{equation}
\label{eq10}
F\left( {z} \right): \equiv \left( {\int_{q} \varphi }  \right)\left( {z}
\right): = \left( {1 - q} \right)z\sum\limits_{k = 0}^{\infty}  {\varphi
\left( {q^{k}z} \right)q^{k}}
\end{equation}
i.e.
\begin{multline}
\label{eq11}
F\left( {z} \right) \equiv \left( {\int_{q} \varphi }
\right)\left( {z} \right) = \left( {1 - q} \right)z\left( {\sum\limits_{k =
0}^{\infty}  {q^{k}\hat {Q}^{k}\varphi} }  \right)\left( {z} \right) =\\
=\left( {\left( {1 - q} \right)z\frac{{1}}{{1 - q\hat {Q}}}\varphi}
\right)\left( {z} \right).
\end{multline}
Of course
\begin{equation}
\label{eq12}
\partial _{q} \circ \int_{q} = id
\end{equation}
as
\begin{equation}
\label{eq13}
\frac{{1 - q\hat{Q}}}{{\left( {1 - q} \right)}}\partial _{0} \left(
{\left( {1 - q}\right)\hat {z}\frac{{1}}{{1 - q\hat {Q}}}} \right)=id.
\end{equation}
\vspace{0.2cm}
\noindent Naturally (\ref{eq13}) might serve to define a right inverse operation to
``$q$-difference-ization''
\begin{center}
 $\left( {\partial _{q} \varphi}  \right)\left( {x} \right) = \frac{{1 - q\hat
{Q}}}{{\left( {1 - q} \right)}}\partial _{0} \varphi \left( {x} \right)$
\end{center}

\noindent and consequently the ``$q$-integration `` as represented by (\ref{eq10}) and
(\ref{eq11}). As it is well known the definite $q$-integral is an numerical
approximation of the definite integral obtained in the $q \to 1$ limit.
Following the $q$-case example we introduce now an $R$-integration
(consult Remark \ref{remtwoone}).
\begin{equation}
\label{eq14}
\int_{R} {x^{n}} = \left( {\hat {x}\frac{{1}}{{R\left( {q\hat {Q}}
\right)}}} \right) x^{n}=
\frac{{1}}{{R\left( {q^{n + 1}} \right)}}x^{n + 1};\quad n \ge 0
\end{equation}
Of course $\partial _{R} \circ \int_{R}= id$ as
\begin{equation}
\label{eq15}
R\left( {q\hat {Q}} \right)\partial _{o}
\left( {\hat {x}\frac{{1}}{{R\left( {q\hat {Q}} \right)}}} \right)=
id .
\end{equation}
\vspace{0.2cm}
\noindent Let us then finally introduce the analogous representation for $\partial
_{\psi}  $ difference-ization
\begin{equation}
\label{eq16}
\partial _{\psi}  = \hat {n}_{\psi} \partial _{o};\; \hat {n}_{\psi}
x^{n - 1}= n_{\psi} x^{n - 1};\; n \ge 1.
\end{equation}

\noindent Then

\begin{equation}
\label{eq17}
\int_{\psi} {x^{n}} = \left( {\hat {x}\frac{{1}}{{\hat
{n}_{\psi} } }} \right)x^{n}= \frac{{1}}{{\left( {n + 1} \right)_{\psi}
}}x^{n + 1};\;n \ge 0
\end{equation}

\noindent and of course
\begin{equation}
\label{eq18}
\partial _{\psi}  \circ \int_{\psi} = id
\end{equation}
{\bf Closing Remark:}

\vspace{0.2cm}
\noindent The picture that emerges discloses the fact that any $\psi
$-representation of finite operator calculus or equivalently - any
$\psi $-representation of GHW algebra makes up an example of the
algebraization of the analysis - naturally when constrained to the
algebra of polynomials. We did restricted all our considerations
to the algebra $P$ of polynomials. Therefore the distinction
in-between difference and differentiation operators disappears.
All linear operators on $P$ are both difference and
differentiation operators if the degree of differentiation or
difference operator is unlimited. For example $\frac{{d}}{{dx}} =
\sum\limits_{k \ge 1} {\frac{{d_{k}} }{{k!}}\Delta ^{k}} $ where
$d_{k} = \left[ {\frac{{d}}{{dx}}x^{\underline {k}} } \right]_{x =
0}  = \left( { - 1} \right)^{k - 1}\left( {k - 1} \right)!$ or
$\Delta =\sum\limits_{n \ge 1} {\frac{{\delta _{n}} }{{n!}}}
\frac{{d^{n}}}{{dx^{n}}}$ where $\delta _{n} = \left[ {\Delta
x^{n}} \right]_{x = 0} =1$. Thus the difference and differential
operators and equations are treated on the same footing. For new
applications - due to the author see \cite{2,3}. Our goal here was
to deliver the general scheme of "$\psi $-umbral" algebraization
of the analysis of general differential operators \cite{11}.
 
\vspace{0.2cm}
\noindent Most of the
general features presented here are known to be pertinent to the
{\it Q} representation of finite operator calculus (Viskov ,
Markowsky, Roman)  where {\it Q}  is any linear operator lowering
degree of any polynomial by one . So it is most general example of
the algebraization of the analysis for general differential
operators \cite{11} .

\section{Glossary}

\noindent Here now come short indicatory glossaries of terms and notation
used by Ward \cite{1}, Viskov \cite{6,7}, Markowsky \cite{11},
Roman \cite{27}-\cite{31} on one side and the Rota-oriented
notation on the other side. See also \cite{32}.
\begin{center}
\begin{longtable}{|c|c|}
          \hline
& \\*
{\bf Ward} & {\bf Rota - oriented} (this note)\\*
& \\*
\endhead
          \hline
& \\*
$\left[n\right];\;\left[n\right]!$ & $n_{\psi};\;n_{\psi}!$\\*
& \\*
basic binomial coefficient $\left[n,r\right]=
\frac{\left[n\right]!}{\left[r\right]!\left[n-r\right]!}$ &
$\psi$-binomial coefficient $\binom{n}{k}_{\psi} \equiv
\frac{n_\psi ^{\underline{k}}}{k_\psi !}$\\*
& \\*
          \hline
& \\*
$D = D_{x}$ - the operator $D$ & $\partial_{\psi}$ - the
$\psi$-derivative\\*
& \\*
$D\,x^n = \left[n\right]\,x^{n-1}$ &
$\partial_{\psi}\,x^n=n_{\psi}\,x^{n-1}$\\*
& \\*
          \hline
&   \\*
$(x+y)^n$ & $(x+_{\psi}y)^n$\\*
&  \\*
$(x+y)^n \equiv \sum\limits_{r=0}^{n}\left[n,r\right]x^{n-r}y^r$  &
$(x+_{\psi} y)^{n}=\sum \limits_{k = 0}^{n}\binom{n}{k}_{\psi}x^ky^{n-k}
$\\*
&  \\*
          \hline

&  \\* basic displacement symbol & generalized shift operator\\* &
\\* $E^t;\;t\in${\bf Z} & $E^{y}\left( {\partial _{\psi} }
\right) \equiv exp_{\psi}  \{y\partial _{\psi} \}$;\;$y\in${\bf
F}\\* &  \\* $E\varphi (x) = \varphi (x+1)$ & $E(\partial

_{\psi})\varphi(x)= \varphi(x+_{\psi}1)$\\* &  \\* $E^t\varphi (x)
= \varphi \left(x+\overline{t}\right)$ & $E^y(\partial _{\psi})x^n
\equiv (x+_{\psi}y)^n$\\* &  \\*
          \hline
& \\*
basic difference operator & $\psi $-difference delta operator\\*
&  \\*
$\Delta = E - id$  & $\Delta _{\psi} = E^y(\partial_{\psi}) - id$\\*
& \\*
$\Delta = \varepsilon (D) - id = \sum\limits_{n=0}^{\infty}
\frac{D^n}{\left[n\right]!} - id$ & \\*
& \\*
          \hline
\end{longtable}
\end{center}
\begin{center}
\begin{longtable}{|c|c|}
          \hline
& \\*
{\bf Roman} & {\bf Rota - oriented} (this note)\\
& \\*
\endhead
          \hline
& \\*
$t;\;tx^n = nx^{n-1}$ & $\partial_{\psi}$ - the $\psi$-derivative\\*
& \\*
& $\partial_{\psi}x^n = n_{\psi}x^{n-1}$\\*
& \\*
$\langle t^{k}|p(x)\rangle = p^{(k)}(0)$ &
$[\partial_{\psi}^kp(x)]|_{x=0}$\\*
& \\*
          \hline
& \\*
evaluation functional & generalized shift operator\\*
& \\*
$\epsilon_{y}(t) = \exp{\{yt\}}$
& $E^y(\partial_{\psi}) = \exp_{\psi}{\{y\partial_{\psi}\}}$\\*
& \\*
$\langle t^k|x^n\rangle = n!\delta_{n,k}$ & \\*
& \\*
$\langle \epsilon_y(t)|p(x)\rangle = p(y)$ &
$[E^y(\partial_{\psi})p_n(x)]|_{x=0} = p_n(y)$\\*
& \\*
$\epsilon_y(t)x^n = \sum\limits_{k \geq 0}\binom{n}{k}
x^ky^{n-k}$ &
$E^y(\partial_{\psi})p_n(x) = \sum \limits_{k \geq 0}\binom{n}{k}_{\psi}
p_k(x)p_{n-k}(y)$\\*
& \\*
          \hline
& \\*
formal derivative & Pincherle derivative \\*
& \\*
$f'(t) \equiv \frac{d}{dt}f(t)$ &
$[Q(\partial_{\psi})]${\bf `}$\equiv
\frac{d}{d\partial_{\psi}}Q(\partial_{\psi})$\\*
& \\*
${\overline f}(t)$ compositional inverse of &
$Q^{-1}(\partial_{\psi})$ compositional inverse of \\*
& \\*
formal power series $f(t)$ & formal power series $Q(\partial_{\psi})$\\*
& \\*
          \hline
& \\*
$\theta_{t};\;\;\theta_{t}x^n = x^{n+1};\;n\geq 0$ &
${\hat x}_{\psi};\;\;{\hat x}_{\psi}x^n =
\frac{n+1}{(n+1)_{\psi}}x^{n+1};\;n\geq 0$ \\*
& \\*
$\theta_{t}t = {\hat x}D$ & ${\hat x}_{\psi}\partial_{\psi} = {\hat x}D =
{\hat N}$\\*
& \\*
          \hline
& \\*
$\sum \limits_{k \geq 0}\frac{s_k(x)}{k_{\psi}!}t^k =$ &
$\sum \limits_{k \geq 0}\frac{s_k(x)}{k_{\psi}!}z^k =$\\*
& \\*
$[g({\overline f}(z))]^{-1}\exp {\{x{\overline f}(t)\}}$ &
$s(q^{-1}(z))\exp_{\psi}{\{xq^{-1}(z)\}}$\\*
& \\*
$\{s_{n}(x)\}_{n \geq 0}$ - Sheffer sequence & $q(t),\,s(t)$ indicators\\*
& \\*
for $(g(t),f(t))$ & of $Q(\partial_{\psi})$ and $S_{\partial_{\psi}}$\\*
& \\*
          \hline
& \\*
$g(t)\,s_{n}(x) = q_{n}(x)$ - sequence &
$s_n(x) = S_{\partial_{\psi}}^{-1}\,q_n(x)$ - $\partial_{\psi}$ - basic\\*
& \\*
associated for $f(t)$ & sequence of $Q(\partial_{\psi})$ \\*
& \\*
          \hline
& \\*
The expansion theorem: & The First Expansion Theorem\\*
& \\*
$h(t) = \sum \limits_{k=0}^\infty \frac{\langle h(t)|p_{k}(x) \rangle}{k!}
f(t)^k$ & $T = \sum \limits_{n\geq 0} \frac{[T\,p_{n}(z)]|_{z=0}}{n_{\psi}}
Q(\partial_{\psi})^n$\\*
& \\*
$p_n(x)$ - sequence associated for $f(t)$ & $\partial_{\psi}$ - basic
polynomial sequence $\{p_n\}_{0}^{\infty}$\\*
& \\*
          \hline
& \\*
$\exp\{y{\overline {f}}(t)\} = \sum \limits_{k=0}^{\infty}\frac{p_{k}(y)}
{k!}t^k$ & $\exp_{\psi}\{xQ^{-1}(x)\} = \sum \limits_{k\geq0}
\frac{p_{k}(y)}{k!}z^k$ \\*
& \\*
          \hline
& \\*
The Sheffer Identity: & The Sheffer ${\psi}$-Binomial Theorem:\\*
& \\*
$s_{n}(x+y) = \sum \limits_{k=0}^{n} \binom{n}{k}p_{n}(y)s_{n-k}(x)$ &
$s_{n}(x+_{\psi} y) = \sum \limits_{k\geq 0}\binom{n}{k}_{\psi}s_{k}(x)
q_{n-k}(y)$\\*
& \\*
          \hline
\end{longtable}
\end{center}
\begin{center}
\begin{longtable}{|c|c|}
          \hline
& \\*
{\bf Viskov} & {\bf Rota - oriented} (this note)\\*
& \\*
\endhead
          \hline
& \\*
$\theta _{\psi }$ - the $\psi$-derivative
& $\partial_{\psi}$ - the $\psi$-derivative\\*
& \\*
$\theta_{\psi}\,x^n = \frac{\psi_{n-1}}{\psi_{n}}x^{n-1}$ &
$\partial_{\psi}\,x^n=n_{\psi}\,x^{n-1}$\\*
& \\*
          \hline
& \\*
$A_p\; (p = \{p_n \}_{0}^{\infty})$  & $Q$\\*
& \\*
$A_p\,p_n = p_{n-1}$ & $Q\,p_n = n_{\psi}p_{n-1}$\\*
& \\*
          \hline
& \\*
$B_p\; (p = \{p_n \}_{0}^{\infty})$  & ${\hat x}_Q$\\*
& \\*
$B_p\,p_n = (n+1)\,p_{n+1}$ & ${\hat x}_Q\,p_n = \frac{n+1}{(n+1)_{\psi}}
p_{n+1}$\\*
& \\*
          \hline
& \\*

$E_{p}^y\; (p = \{p_n \}_{0}^{\infty})$ & $E^{y}\left( {\partial
_{\psi} } \right) \equiv exp_{\psi}  \{y\partial _{\psi} \}$\\* &
\\* $E_{p}^y\,p_{n}(x) =\sum \limits_{k=0}^{n}p_{n-k}(x)p_{k}(y)$
& $E^{y}\left( {\partial _{\psi} }  \right)\,p_{n}(x) = $\\* & \\*
& $=\sum \limits_{k \geq 0}
\binom{n}{k}_{\psi}\,p_{k}(x)p_{n-k}(y)$\\* & \\*
           \hline
& \\*

$T - \varepsilon_{p}$-operator: & $E^y$ - shift operator:\\* & \\*
$T\,A_{p} = A_{p}\,T$  & $E^y\varphi (x) = \varphi
(x+_{\psi}y)$\\* & \\*
           \hline
& \\*
& $T$ - $\partial_{\psi}$-shift invariant operator:\\*
& \\*
$\forall_{y\in F}\;TE_{p}^y = E_{p}^{y}T$ &
 $\forall_{\alpha \in F}\;[T,E^{\alpha}(\partial_{\psi})]=0$\\*
& \\*
           \hline
& \\*
$Q$ - $\delta_{\psi}$-operator: &
$Q(\partial_{\psi})$ - $\partial_{\psi}$-delta-operator:\\*
& \\*
$Q$ - $\epsilon_{p}$-operator and &
$Q(\partial_{\psi})$ - $\partial_{\psi}$-shift-invariant and\\*
& \\*
$Qx = const \neq 0$ &
$Q(\partial_{\psi})(id) = const \neq 0$\\*
& \\*
          \hline/
& \\*
$\{p_n(x),n\geq 0\}$ - $(Q,\psi)$-basic &
$\{p_n\}_{n\geq 0}$ -$\partial_{\psi}$-basic\\*
& \\*
polynomial sequence of the &
polynomial sequence of the\\*
& \\*
$\delta_{\psi}$-operator $Q$ &
$\partial_{\psi}$-delta-operator $Q(\partial_{\psi})$\\*
& \\*
          \hline
& \\*
$\psi$-binomiality property & $\psi$-binomiality property\\*
& \\*
$\Psi_{y}s_n(x) = $ & $E^y(\partial_{\psi})p_{n}(x) = $\\*
& \\*
$=\sum \limits_{m=0}^n \frac{\psi_{n}\psi_{n-m}}{\psi_{n}}s_m(x)
p_{n-m}(y)$ &
$= \sum \limits_{k \geq 0}\binom{n}{k}_{\psi}p_{k}(x)p_{n-k}(y)$\\*
& \\*
          \hline
& \\*
$T = \sum \limits_{n \geq 0}\psi_n[VTp_n(x)]Q^n$ &
$T = \sum \limits_{n\geq 0} \frac{[T p_{n}(z)]|_{z=0}}{n_{\psi}!}
Q(\partial_{\psi})^n$\\*
& \\*
$T\Psi_{y}p(x) =$ & $T p(x+_{\psi} y) =$\\*
& \\*
$\sum \limits_{n \geq 0}\psi_ns_n(y)Q^nSTp(x)$  &
 $\sum \limits_{k \geq 0} \frac{s_{k}(y)}
{k_{\psi}!} Q(\partial_{\psi})^{k}ST p(x)$ \\*
& \\*
          \hline
\end{longtable}
\end{center}
\begin{center}
\begin{longtable}{|c|c|}
          \hline
& \\*
{\bf Markowsky} & {\bf Rota - oriented}\\*
& \\*
          \endhead
          \hline
& \\*
$L$ - the differential operator & $Q$\\*
& \\*
$L\,p_n = p_{n-1}$ & $Q\,p_n = n_{\psi}p_{n-1}$\\*
& \\*
          \hline
& \\*
$M$ & ${\hat x}_Q$\\*
& \\*
$M\,p_n = p_{n+1}$  & ${\hat x}_Q\,p_n = \frac{n+1}{(n+1)_{\psi}}
p_{n+1}$\\*
& \\*
          \hline
& \\*
$L_y$ &
$E^{y}\left( Q \right)= \sum \limits_{k \geq 0} \frac{p_k (y)}{k_{\psi}!}
Q^k$\\*
& \\*
$L_y\,p_{n}(x) = $ &
$E^{y}\left( Q \right)\,p_{n}(x) = $\\*
& \\*
$= \sum \limits_{k=0}^{n} \binom{n}{k}p_{k}(x)p_{n-k}(y)$ &
$=\sum \limits_{k \geq 0} \binom{n}{k}_{\psi}\,p_{k}(x)p_{n-k}(y)$\\*
& \\*
           \hline
& \\*
$E^a$ - shift-operator: & $E^y$ - $\partial_{\psi}$-shift operator:\\*
& \\*
$E^a\,f(x) = f(x+a)$ & $E^y\varphi (x) = \varphi (x+_{\psi}y)$\\*
& \\*
           \hline
& \\*
$G$ - shift-invariant operator:
& $T$ - $\partial_{\psi}$-shift invariant operator:\\*
& \\*
$EG=GE$ &
 $\forall_{\alpha \in F}\;[T,E(Q)]=0$\\*
& \\*
           \hline
& \\*
$G$ - delta-operator: &
$L=L(Q)$ - $Q_{\psi}$-delta operator:\\*
& \\*
$G$ - shift-invariant and &
$[L,Q] = 0$ and\\*
& \\*
$Gx = const \neq 0$ & $L(id) = const \neq 0$\\*
& \\*
          \hline
& \\*
$D_L(G)$ & $G\mathbf{'} = [G(Q),{\hat x}_{Q}]$ \\*
& \\*
$L$ - Pincherle derivative of $G$ & $Q$ - Pincherle derivative\\*
& \\*
$D_L(G) = [G,M]$ & \\*
& \\*
          \hline
& \\*
$\{Q_{0},Q_{1},...\}$ - basic family &
$\{p_n\}_{n\geq 0}$ -$\psi$-basic\\*
& \\*
for differential operator $L$ & polynomial sequence of the\\*
& \\*
& generalized difference operator $Q$\\*
& \\*
          \hline
& \\*
binomiality property & $Q$ - $\psi$-binomiality property\\*
& \\*
$P_n(x+y) = $  & $E^y(Q)p_{n}(x) = $\\*
& \\*
$=\sum \limits_{i=0}^n \binom{n}{i}P_i(x)P_{n-i}(y)$ &
$= \sum \limits_{k \geq 0}\binom{n}{k}_{\psi}p_{k}(x)p_{n-k}(y)$\\*
& \\*
          \hline
\end{longtable}
\end{center}
\vspace{0.2cm}
\noindent
{\bf {\em Acknowledgements}}

Discussions with Professor Oleg V. Viskov are highly acknowledged.
The author expresses his gratitude  to Ewa Borak for help in preparation
of the LaTeX version of this contribution.
The author is also very much indebted to the Referee whose indications allowed preparing
    the paper in hopefully

     more desirable form.


\begin{thebibliography}{10}

\bibitem{4}
R. P. Boas and Jr. R. C. Buck: {\em Am. Math. Monthly\/} {\bf 63},
626 (1959).

\bibitem{5}
R. P. Boas and Jr. R. C. Buck: {\em Polynomial Expansions of
Analytic Functions\/}, Springer, Berlin 1964.

\bibitem{14}
A. Di Bucchianico and D.Loeb: {\em J. Math. Anal. Appl.\/} {\bf
92}, 1 (1994).

\bibitem{21}
A. Di Bucchianico and D.Loeb: {\em Integral Transforms and Special
Functions\/} {\bf 4}, 49 (1996).

\bibitem{23}
A. Di Bucchianico and D.Loeb: {\em J. Math. Anal. Appl.\/} {\bf
199}, 39 (1996).

\bibitem{17}
P. Feinsilver and R. Schott: {\em Algebraic Structures and
Operator Calculus\/}, Kluwer Academic Publishers, New York 1993.

\bibitem{16}
C. Graves: {\em Proc. Royal Irish Academy\/} {\bf 6}, 144
(1853-1857).


\bibitem{24}
F. H. Jackson: {\em Quart. J. Pure and Appl. Math.\/} {\bf 41},
193 (1910).

\bibitem{25}
F. H. Jackson: {\em Messenger of  Math.\/} {\bf 47}, 57 (1917).

\bibitem{26}
F. H. Jackson: {\em Quart. J. Math.\/} {\bf 2}, 1 (1951).

\bibitem{22}
P. Kirschenhofer: {\em Sitzunber. Abt. II  Oster. Ackad. Wiss.
Math. Naturw. Kl.\/} {\bf 188}, 263 (1979).

\bibitem{20}
S. G. Kurbanov and V. M. Maximov: {\em Dokl. Akad. Nauk Uz.
SSSR\/} {\bf 4}, 8 (1986).

\bibitem{2}
A. K. Kwa\'sniewski \textit{ Rep. Math. Phys.} \textbf{48} (3),
305 (2001)

\bibitem{3}
A. K. Kwa{\'s}niewski \textit{Integral Transforms and Special
Functions}\textbf{2} (4), 333 (2001)

\bibitem{12}
A. K. Kwa\'sniewski: {\em Advances in Applied Clifford Algebras\/}
{\bf 9}, 41 (1999).

\bibitem{32}
A. K. Kwasniewski and E. Gradzka: {\em  Rendiconti del Circolo
Matematico di Palermo Serie II , Suppl.\/} {\bf 69}, 117(2002).

\bibitem{11}
G. Markowsky: {\em J. Math. Anal. Appl.\/} {\bf 63}, 145 (1978).

\bibitem{19}
S. Pincherle and U. Amaldi: {\em Le operazioni distributive e le
loro applicazioni all`analisi\/},  N. Zanichelli, Bologna 1901.

\bibitem{27}
S. M. Roman: {\em J. Math. Anal. Appl.\/} {\bf 87}, 58 (1982).

\bibitem{28}
S. M. Roman: {\em J. Math. Anal. Appl.\/} {\bf 89}, 290 (1982).

\bibitem{29}
S.M. Roman: {\em J. Math. Anal. Appl.\/} {\bf 95}, 528 (1983).

\bibitem{30}
S. M. Roman: {\em The umbral calculus\/}, Academic Press, New York
1984.

\bibitem{31}
S. R. Roman: {\em J. Math. Anal. Appl.\/} {\bf 107}, 222 (1985).

\bibitem{8}
G.-C. Rota and R. Mullin: {\em On the foundations of combinatorial
theory, III. Theory of Binomial  Enumeration in "Graph Theory and
Its Applications"\/}, Academic Press, New York 1970.

\bibitem{9}
G. C. Rota, D.Kahaner and A. Odlyzko: {\em J. Math. Anal. Appl.\/}
{\bf 42}, 684 (1973).

\bibitem{10}
G. C. Rota: {\em Finite Operator Calculus\/}, Academic Press, New
York 1975.

\bibitem{15}
N. Ya. Sonin: {\em Izw. Akad. Nauk\/} {\bf 7}, 337 (1897).

\bibitem{33}
J.F. Steffenson {\em Acta Mathematica\/} {\bf 73}, 333 (1944).

\bibitem{6}
O.V. Viskov: {\em Soviet Math. Dokl.\/} {\bf 16}, 1521 (1975).

\bibitem{7}
O.V. Viskov: {\em Soviet Math. Dokl.\/} {\bf 19}, 250 (1978).

\bibitem{13}
O.V. Viskov: {\em Trudy Matiematicz`eskovo Instituta AN SSSR\/} {\bf 177},
21 (1986).

\bibitem{18}
O.V. Viskov: {\em Integral Transforms and Special Functions\/} {\bf 1},
2 (1997).

\bibitem{1}
M. Ward: {\em Amer. J. Math.\/} {\bf 58}, 255 (1936).


\end{thebibliography}
\end{document}